\documentclass{article}
\usepackage{amsmath,amssymb,amsfonts}
\usepackage{algorithmic}
\usepackage{graphicx}
\usepackage[numbers]{natbib}
\usepackage{graphicx,listings}
\usepackage{textcomp}
\usepackage{xcolor}
\usepackage{tcolorbox}
\usepackage{multicol}
\usepackage{tcolorbox}
\usepackage{scrextend}
\usepackage{alltt}
\usepackage{pifont}
\usepackage{mathtools,array}
\usepackage{tikz}
\usepackage{tkz-euclide}
\usepackage{etoolbox}
\usepackage{pgfplots}
\usepackage{pifont}
\usetikzlibrary{calc,angles,positioning,intersections,quotes,decorations.markings, through,backgrounds,turtle}
\usepackage{tikz-3dplot}
\def\disp{\displaystyle}
\newcommand{\mysim}{\mbox{$\asymp$}}
\def\myfrac#1#2{\frac{\disp{#1}}{\disp{#2}}}

\def\BibTeX{{\rm B\kern-.05em{\sc i\kern-.025em b}\kern-.08em
    T\kern-.1667em\lower.7ex\hbox{E}\kern-.125emX}}

\newtheorem{example}{Example}
\newtheorem{theorem}{Theorem}
 
 \newtheorem{definition}{Definition}
 
\begin{document}

\title{The paradoxes and the infinite dazzled ancient
 mathematics and continue to do so today}

  

\author{Fairouz Kamareddine and Jonathan Seldin}
\maketitle

\begin{abstract}
  This paper looks at how ancient mathematicians (and especially the Pythagorean school) were faced by problems/paradoxes associated with  the infinite which led them to juggle two systems of numbers: the discrete whole/rationals which were handled arithmetically and the continuous magnitude quantities which were handled geometrically.  We look at how  approximations and mixed numbers (whole numbers with fractions) helped develop the arithmetization of geoemtry and  the development of mathematical analysis and real numbers.
  \end{abstract}


\section{Why did it take so long to develop Real Numbers and Analysis?}
 \begin{center}   \begin{tcolorbox}[width=.88\textwidth, colframe=red]
     God made the integers; all else is the work of man.
  \mbox{} $\hspace{1in} \phantom \quad$ \hfill  Kronecker
 \end{tcolorbox}\end{center}
 The concepts and language of mathematics have been under development slowly but surely since ancient times.
 Despite the obstacles, this development
 uncovered  fascinating results, which include  as late as the 20th century, a sound foundation of the theory of the infinitesimal (which is in essence the foundation of mathematics) and the theory of the computable. Well before then,  Leibniz (1646--1717) conceived of {\em  automated deduction} where he wanted to find  a language $L$
 and a method  that could carry out proof checking/finding to determine the correctness of statements in $L$.\footnote{Now we know, due to
  later results by G{\"o}del, Church and Turing, that such a
  method can not work for every statement.}  Leibniz was frustrated by the limitations in expressing thoughts:
 \begin{center}   \begin{tcolorbox}[width=.88\textwidth, colframe=red]
  If we could find characters or signs appropriate for expressing all our thoughts as definitely and as exactly as arithmetic expresses numbers or geometric analysis expresses lines, we could in all subjects in so far as they are amenable to reasoning accomplish what is done in  Arithmetic/Geometry.
  \mbox{} $\hspace{1in} \phantom \quad$ \hfill  Leibniz
 \end{tcolorbox}\end{center}
 But at the time of Leibniz, expressibility in Arithmetic was far from complete and the real numbers were still not developed.
 The later development of real analysis\footnote{Thanks to Euler who converted the calculus of  Newton and Leibniz from a geometrical field to a field where mathematical formulae are analysed.} would be based on the real numbers and 
 the arithmetisation of geometry.
 
\subsection{From naturals to intergers and rationals}
        {\em Natural numbers}
                were long understood,  but it may come as a surprise that as late as the 14th century, negative numbers were not known in Europe.  
              In Italy, a double entry bookkeeping system compensated for their absence.
        Accounts in which debits may be greater than  credits were compared without using negative integers.
 If  $c$ and $d$ are in $\mathbb{N}^+$, then \emph{account} $c \ominus d$ 
 has credit $c$ and debit $d$. Define $\mbox{accounts} = \{m \ominus n\: | \: m,n,p,q\in \mathbb{N}^+\}$.  Just like the arithmetic $(\mathbb{N}^+, =, +, \cdot, 1)$ on natural numbers $\mathbb{N}^+=\{1,2,\cdots\}$ 
 is defined with equality =, addition $+$, multiplication $\cdot$, and identity element 1 for $\cdot$, we define  $(\mbox{accounts}, \cong, +_c, \cdot_c)$ by:
  \begin{itemize}
\item
    $m \ominus n \cong p \ominus q \mbox{ iff } m + q = n + p .$
  \item
    $(m \ominus n) +_c (p \ominus q) = (m + p) \ominus (n + q)$.
    \item
    $(m \ominus n)\cdot_c (p \ominus q) = (mp +  nq) \ominus (mq + np) .$
  \end{itemize}
  The integers $(\mathbb{Z},  +_i, \cdot_i, 0_i, 1_i, -\alpha)$ are then defined from the equivalence classes: $\left[ m \ominus n \right] = \left\{ p \ominus q : p \ominus q \cong m \ominus n \right \}$ by:
  \begin{itemize}
  \item
    $\mathbb{Z} = \{[m\ominus n] \: | \:  m, n \in \mathbb{N}^+\}$.
    \item
      $ [(m \ominus n)] +_i [(p \ominus q)] = [(m \ominus n) +_c (p \ominus q)]$.
    \item
      $[(m \ominus n)] \cdot_i [(p \ominus q)] =[(m \ominus n) \cdot_c (p \ominus q)]$.
\item
  Identity $0_i$ for $+_i$: for any $m$, $n$ in $\mathbb{N}^+$, $[m \ominus m] = [n \ominus n]$.
  \item
    Identity $1_i$ for $\cdot_i$: take  $1_i = [(p+1) \ominus p]$ for any $p\in \mathbb{N}^+$.
    \item
   Inverse for $+_i$: if $\alpha = [m \ominus n]$, then $- \alpha = [n \ominus m]$.
  \end{itemize}
  Like we
  defined $(\mathbb{Z},  +_i, \cdot_i, 0_i, 1_i, -\alpha)$ from $(\mbox{accounts}, \cong, +_c, \cdot_c)$ which were defined from $(\mathbb{N}^+, =, +, \cdot, 1)$,
  we can define positive rational numbers $(\mathbb{Q}^+,  +_r, \cdot_r, 1_r, \mathbf{a}^{-1})$ from $\mbox{fractions}=\{\myfrac{m}{n} \: | \: m,n\in \mathbb{N}^+\}$ where the arithmetic of 
      $(\mbox{fractions}, \mysim, +_f, \cdot_f)$ is defined from $(\mathbb{N}^+, =, +, \cdot, 1)$ by:
\begin{itemize}
       \item
         $\myfrac{m}{n} \mysim \myfrac{p}{q} \mbox{ if and only if }  mq = np$,
         \item
           $\myfrac{m}{n} +_f \myfrac{p}{q} = \myfrac{mq +np}{nq}$ and  $\myfrac{m}{n} \cdot_f \myfrac{p}{q} = \myfrac{mp}{nq}$.
             \end{itemize}    
Then,         we     define $(\mathbb{Q}^+,  +_r, \cdot_r, 1_r, \mathbf{a}^{-1})$ from 
equivalence classes $\left[ \myfrac{m}{n} \right] = \left\{ \myfrac{p}{q} \: | \:  \myfrac{p}{q} \mysim \myfrac{m}{n}  \right\}$ as follows:

 \begin{itemize}
               \item
                      $\mathbb{Q}^+ = \{\left[ \myfrac{m}{n} \right] \: | \: m,n \in \mathbb{N}^+\}$.
          \item
            $\left[\myfrac{m}{n}\right] +_r \left[\myfrac{p}{q}\right] = \left[\myfrac{m}{n} +_f \myfrac{p}{q} \right]$
            and 
           $\left[\myfrac{m}{n}\right] \cdot_r \left[\myfrac{p}{q}\right] = 
       \left[ \myfrac{m}{n} \cdot_f \myfrac{p}{q} \right]$.
     \item
       Identity $1_r$ for $\cdot_r$: take  $1_r = \left[ \myfrac{1}{1} \right]$.
       \item
   Inverse for $+_r$: $\left[ \myfrac{m}{n} \right]^{-1} = \left[ \myfrac{n}{m} \right]$.
 \end{itemize}
 The steps to build $(\mathbb{Z},  +_i, \cdot_i, 0_i, 1_i, -\alpha)$ and $(\mathbb{Q}^+,  +_r, \cdot_r, 1_r, \mathbf{a}^{-1})$  lead to a generalisation as follows:
 \begin{definition}
     If  $(S,\circ)$ is a  Commutative Cancellation Semigroup\footnote{I.e.,  $\circ$  satisfies closure, commutativity, associativity and cancellation law on $S$ where cancellation means that $a\circ b = a\circ c$ implies $b = c$.} (CCS), then build $(S\times S, \approx, *)$ as follows:
          \begin{itemize}
          \item
           Define congruence $\approx$ on $S\times S$ based on $(S, \circ)$ by:
            $(x,y) \approx (u,v)$ iff $x \circ v = y \circ u$.
                     \item
                                             The operation $*$ on $S\times S$ inherited from $\circ$ is defined by                        
                                             $(x,y) * (u,v) = (x \circ u,y \circ v)$.
          \end{itemize}
          Then,                       define $[(x,y)] = \{ (u,v) : (u,v) \approx (x,y) \}$ and   $S_d = \{[(x,y)] : x, y\in S\}$, and build $(s_d, \circ_d, e_d, \mathfrak{a}^{-1})$ as  follows:

                                                     \begin{itemize}
                                          \item
                                            Define
$[(x,y)]e_d[(u,v)]=[(x,y)*(u,v)]= [(x \circ u, y \circ v)]$.  Note that 
  $(S_d, \circ_d)$ is a CCS.
\item Note that if $x\in S$, then $x_d= [(y \circ x, y)] \in S_d$.
 \item  Identity:
  Define $e_d$ to be $[(x,x)]$ for some $x$ in $S$.  For all  $\mathfrak{a}$, we have
$e_d \circ_d \mathfrak{a} = \mathfrak{a} \circ_d e_d = \mathfrak{a}$.
  \item  Inverses:
  If  $\mathfrak{a} = [(x,y)]$, define  $\mathfrak{a} ^{-1}$ to be $[(y,x)]$.
    We have
    $ \mathfrak{a} \circ_d \mathfrak{a}^{-1} = e_d = \mathfrak{a}^{-1} \circ_d \mathfrak{a}$.
\end{itemize}
 \end{definition}
  Comparing  the theory of fractions and the theory of accounts suggests that we can define a unified theory for adding inverses and, if none is present, identity elements. 
\begin{center}
 \begin{tabular}{|l|l|l|}
  \hline
  CCS &$(\mathbb{N}^+, +)$ & $(\mathbb{N}^+,\cdot)$ \\
  \hline
  inverses &$\times$ & $\times$  \\
  \hline
  Identity  element &$\times$ & $\surd$  \\
  \hline
  CCS& $(\mathbb{Z}, +_i)$ & $(\mathbb{Q}^+,\cdot_r)$ \\
  with identity and inverses &$\surd$ & $\surd$\\
  \hline
\end{tabular}
\end{center}
Just like we built $(\mathbb{Z}, +)$ with identity $0_i$ and inverses $-\alpha$ from $(\mathbb{N}^+, +)$,  we can build $(\mathbb{Q}, +_{r_i})$ with identity and inverses from $(\mathbb{Q}^+, +_r)$.  But we cannot build $\mathbb{R}$ this way.  The real numbers need to be constructed (using approximations and limits like Dedekind cuts, Cauchy sequences, etc.).   This brings us to what is the foundations of mathematics?  The foundation of mathematics is reasoning about whether the infinitesimal is sound.    Euclid’s Elements developed mathematics in geometric terms and anything not expressible in such terms was excluded.   Geometry could accommodate the whole numbers and their ratios as well as irrational magnitudes. Think for example of  the spiral of Theodorus of Cyrene which established that the square roots of non square integers from 3 to 17 are irrationals.
\begin{center}
  \begin{tikzpicture}[scale=0.75]
  \def\numOfTs{15}
  \def\dotRad{1pt}
  \def\angleMarkLen{4pt}
  \draw [turtle={home,rt,fd,lt,fd}] coordinate (p0)
  \foreach \i in {1,...,\numOfTs}
  {%
    [turtle={%
      lt={90-atan(sqrt(\i))},%
      fd=\angleMarkLen,lt,%
      fd=\angleMarkLen,lt,%
      fd=\angleMarkLen,lt,%
      fd=\angleMarkLen,lt,fd}] coordinate (p\i)%
  };
  \foreach[count=\j from 2] \i in {0,1,...,\numOfTs}
  {%
    \draw (0,0) -- node[fill=white]{\tiny$\sqrt{\j}$} (p\i);%
    \fill (p\i) circle (\dotRad);%
  }
  \fill (1,0) circle (\dotRad);
  \fill (0,0) circle (\dotRad);
  \end{tikzpicture}
  \end{center}
\subsection{Proofs by Pebbles/Diagrams}
    Knorr~\cite{knorr:EEEl} suggests that the original proofs were proofs as diagrams using \emph{pebble diagrams}.
    It is known that the ancient Greeks did arithmetic by counting with pebbles, and pebble diagrams give these calculations by representing the  pebbles by using small circles.
    \begin{example}
      Here are some statements and their  proofs:
      \begin{itemize}
        \item
      The square of an odd number is 1 +  a multiple of 4.\\
      The square of an even number is a multiple of 4.
      \begin{center}
      \begin{tikzpicture}[scale=0.4]
    \filldraw
    (5.5,2.5) -- (2.5,2.5);
        \filldraw
    (0.5,1.5) -- (3.5,1.5);
     \filldraw
   (3.5,-0.5) -- (3.5,2.5);
      \filldraw
      (2.5,1.5) -- (2.5,4);
      \draw (1,0) circle (.1);
      \draw (2,0) circle (.1);
      \draw (3,0) circle (.1);
      \draw (4,0) circle (.1);
      \draw (5,0) circle (.1);
      \draw (1,1) circle (.1);
      \draw (2,1) circle (.1);
      \draw (3,1) circle (.1);
      \draw (4,1) circle (.1);
      \draw (5,1) circle (.1);

      \draw (1,2) circle (.1);
      \draw (2,2) circle (.1);
      \draw (3,2) circle (.1);
      \draw (4,2) circle (.1);
      \draw (5,2) circle (.1);

      \draw (1,3) circle (.1);
      \draw (2,3) circle (.1);
      \draw (3,3) circle (.1);
      \draw (4,3) circle (.1);
      \draw (5,3) circle (.1);
         
      \draw (1,4) circle (.1);
      \draw (2,4) circle (.1);
      \draw (3,4) circle (.1);
      \draw (4,4) circle (.1);
      \draw (5,4) circle (.1);
      \end{tikzpicture}
      \hspace{0.2in}
     \begin{tikzpicture}[scale=0.4]
    \filldraw
    (3.5,-0.5) -- (3.5,5.5);
      \filldraw
      (0.5,2.5) -- (6.5,2.5);
      \draw (1,0) circle (.1);
      \draw (2,0) circle (.1);
      \draw (3,0) circle (.1);
      \draw (1,1) circle (.1);
      \draw (2,1) circle (.1);
      \draw (3,1) circle (.1);
      \draw (1,2) circle (.1);
      \draw (2,2) circle (.1);
      \draw (3,2) circle (.1);
      
      \draw (4,0) circle (.1);
      \draw (5,0) circle (.1);
      \draw (6,0) circle (.1);
      \draw (4,1) circle (.1);
      \draw (5,1) circle (.1);
      \draw (6,1) circle (.1);
      \draw (4,2) circle (.1);
      \draw (5,2) circle (.1);
      \draw (6,2) circle (.1);

      \draw (1,3) circle (.1);
      \draw (2,3) circle (.1);
      \draw (3,3) circle (.1);
      \draw (1,4) circle (.1);
      \draw (2,4) circle (.1);
      \draw (3,4) circle (.1);
      \draw (1,5) circle (.1);
      \draw (2,5) circle (.1);
      \draw (3,5) circle (.1);
      
       \draw (4,3) circle (.1);
      \draw (5,3) circle (.1);
      \draw (6,3) circle (.1);
      \draw (4,4) circle (.1);
      \draw (5,4) circle (.1);
      \draw (6,4) circle (.1);
      \draw (4,5) circle (.1);
      \draw (5,5) circle (.1);
      \draw (6,5) circle (.1);
     \end{tikzpicture}
     \end{center}
  \item
    \emph{If as many odd numbers as we please be
added together, and their multitude be even, then the
sum is even.}
      \begin{center}
  \begin{tikzpicture}[scale=0.2]
    \path (0,0) coordinate (A) (1,0) coordinate (A1)(2,0) coordinate (A2)(3,0) coordinate (A3)(4,0) coordinate (A4)(5,0) coordinate (A5)
    (6,0) coordinate (B)(7,0) coordinate (B1)
    (8,0) coordinate (B2) (9,0) coordinate (C1) (10,0) coordinate (C2)(11,0) coordinate (C3)(12,0) coordinate (C4)(13,0) coordinate (C5)(14,0) coordinate (C6)(15,0) coordinate (C7)
    (16,0) coordinate (C)(17,0) coordinate (D0)(18,0) coordinate (D)(19,0) coordinate (D1) (20,0) coordinate (E1)(21,0) coordinate (E2)(22,0) coordinate (E3)(23,0) coordinate (E4)(24,0) coordinate (E5)(25,0) coordinate (E6)
    (26,0) coordinate (E);
        \draw (A)
    -- (B) node [at start, above ] {$A$}
    -- (C) node [at start, above ] {$B$}
    -- (D) node [at start, above ] {$C$}
    -- (E) node [at start, above ] {$D$}
    -- (A) node [at start, above ] {$E$}
    -- cycle;
    \tkzDrawPoints[color=red](A,B,C,D,E);
      \end{tikzpicture}      
    \end{center}

      \vspace{0.1in}
       \begin{center}
  \begin{tikzpicture}[scale = 0.5]
\filldraw 
(0,7) circle (2pt) --
(1,7) circle (2pt) --
(2,7) circle (2pt) --
(3,7) circle (2pt) --
(4,7) circle (2pt);
\filldraw 
(6,7) circle (2pt) --
(7,7) circle (2pt) --
(8,7) circle (2pt) --
(9,7) circle (2pt) --
(10,7) circle (2pt);
\filldraw
(5,7) circle (2pt) --
(5,8) circle (2pt);
\filldraw 
(2,8) circle (2pt) --
(3,8) circle (2pt) --
(4,8) circle (2pt);
\filldraw 
(6,8) circle (2pt) --
(7,8) circle (2pt) --
(8,8) circle (2pt);
\filldraw
(1,5) circle (2pt) --
(2,5) circle (2pt) --
(3,5) circle (2pt) --
(4,5) circle (2pt);
\filldraw 
(6,5) circle (2pt) --
(7,5) circle (2pt) --
(8,5) circle (2pt) --
(9,5) circle (2pt);
\filldraw
(5,5) circle (2pt) --
(5,6) circle (2pt);
\filldraw 
(4,6) circle (2pt);
\filldraw 
(6,6) circle (2pt);

\filldraw 
(0,2) circle (2pt) --
(1,2) circle (2pt) --
(2,2) circle (2pt) --
(3,2) circle (2pt) --
(4,2) circle (2pt);
\filldraw 
(6,2) circle (2pt) --
(7,2) circle (2pt) --
(8,2) circle (2pt) --
(9,2) circle (2pt) --
(10,2) circle (2pt);
\filldraw 
(2,3) circle (2pt) --
(3,3) circle (2pt) --
(4,3) circle (2pt);
\filldraw 
(6,3) circle (2pt) --
(7,3) circle (2pt) --
(8,3) circle (2pt);
\filldraw
(1,0) circle (2pt) --
(2,0) circle (2pt) --
(3,0) circle (2pt) --
(4,0) circle (2pt);
\filldraw 
(6,0) circle (2pt) --
(7,0) circle (2pt) --
(8,0) circle (2pt) --
(9,0) circle (2pt);
\filldraw 
(4,1) circle (2pt);
\filldraw 
(6,1) circle (2pt);
\filldraw 
(6,-1) circle (2pt) --
(7,-1) circle (2pt);
\filldraw 
(4,-1) circle (2pt) --
(3,-1) circle (2pt);
  \end{tikzpicture}
  \end{center}
   \end{itemize}
    \end{example}
    The Greeks also mastered the use of geometric proofs:
    \begin{example}
      The geometric proof of the Pythagorean Theorem:  $c^2 = a^2 +b^2$.
      The left square shows $(a+b)^2 = 2ab+c^2$ 
      while the right one shows $(a+b)^2 = 2ab+a^2+b^2$.\\  Hence, $2ab+c^2 = 2ab+a^2+b^2$ and  $c^2 = a^2 +b^2$.
        \begin{center}
 \begin{tikzpicture}[scale = 0.55]
    \path (0,3) coordinate (H) (0,0) coordinate (A) (1,0) coordinate (E);
    \draw (A)
      -- (H) node [midway, left] {$b$}
      -- (E) node [pos=0.5, right] {$c$}
      -- (A) node [midway, below]  {$a$}
      -- cycle;
       \path  (4,0) coordinate (B) (4,1) coordinate (F);
    \draw (E)
      -- (F) node [midway, above] {$c$}
      -- (B) node [midway, right] {$a$}
      -- (E) node [midway, below]  {$b$}
      -- cycle;
           \path  (4,4) coordinate (C) (3,4) coordinate (G);
    \draw (F)
      -- (C) node [midway, right] {$b$}
      -- (G) node [midway, above] {$a$}
      -- (F) node [midway, left]  {$c$}
      -- cycle;
           \path  (0,4) coordinate (D);
    \draw (D)
      -- (G) node [midway, above] {$b$}
      -- (H) node [midway, below] {$c$}
      -- (D) node [midway, left]  {$a$}
      -- cycle;
         \path (6,4) coordinate (L) (9,4) coordinate (O) (6,3) coordinate (P);
    \draw (L)
      -- (O) node [midway, above] {$b$}
      -- (P) 
      -- (L) node [midway, left]  {$a$}
      -- cycle;
       \path  (10,4) coordinate (K) (10,3) coordinate (N) (9,3) coordinate (Z);
    \draw (O)
      -- (K) node [midway, above] {$a$}
      -- (N) node [midway, right] {$a$}
      -- (Z)
      -- (O)
      -- cycle;
           \path  (10,0) coordinate (J) (9,0) coordinate (M);
    \draw (N)
      -- (J) node [midway, right] {$b$}
      -- (M) node [midway, below] {$a$}
      -- (N)
      -- cycle;
              \path  (6,0) coordinate (I);
        \draw (P)
      -- (I) node [midway, left] {$b$}
      -- (M) node [midway, below] {$b$}
      -- (Z)
      -- (P)
      -- cycle;
 \end{tikzpicture}
      \end{center}
    \end{example}
    \subsection{Proofs by Contradiction}
    According to Knorr~\cite{knorr:EEEl},    the change from  proofs using diagrams/pebbles to proofs as sequences of statements occurred with  the  discovery of incommensurability:
\begin{tcolorbox}[width=.88\textwidth, colframe=red]
  \begin{theorem}
    \label{incommensthe}
 There is no unit which measures exactly the side and diagonal of a square.
    \end{theorem}
\end{tcolorbox} 
Key results needed for the incommensurability proof relate to Pythagorean triples and the theory of Odd/Even Numbers:
      \begin{definition}
      {\it Pythagorean triples}
    are triples of positive whole numbers representing the lengths of two legs and the hypotenuse of a right triangle. 
    I.e., a Pythagorean triple is a triple of positive integers $(a, b, c)$ if and only if $a^2 + b^2 = c^2$.
   \end{definition}
      E.g.\ (3, 4, 5), (6, 8, 10), (5, 12, 13), (9, 12, 15), (8, 15, 17).
      
      The following are  the results needed to prove incommensurability  theorem~\ref{incommensthe}.
      Assume $(a, b, c)$ is a Pythagorean triple. 
      
  \begin{itemize}
  \item[1.]
    \label{itemone}
      If $c$ is even, then both $a$
      and $b$ are even.
    \item[2.]
    \label{itemtwo}
    If $c$ is even, then $(\myfrac{a}{2} , \myfrac{b}{2} , \myfrac{c}{2} )$
      is also a Pythagorean triple.
    \item[3.]
    \label{itemthree}
    If $c$ is a multiple of four, then 
      so are $a$ and $b$.
    \item[4.]
    \label{itemfour}
    If $c$ is odd, then one of $a$, $b$ is odd and the other is even.
    \item[5.]
    \label{itemfive}
    If any two of $a, b, c$ 
      is even, then the third is also even.
    \item[6.]
    \label{itemsix}
    If one of $a, b, c$ 
    is odd, then two 
    are odd and one is even. 
  \end{itemize}
   $1\cdots 6$ above can be shown using diagrams/pebbles.
However, theorem~\ref{incommensthe} itself needs a proof by contradiction:
   
Proof. Suppose there is such a unit in terms of which, the side of the square is $a$ and the diagonal is $c$. \\
    Then, we have a right triangle \begin{tikzpicture}[scale=0.25]
      \coordinate  (A) at (-1cm,-1.cm);
\coordinate  (C) at (1.5cm,-1.0cm);
\coordinate  (B) at (1.5cm,1.0cm);
\draw (A) -- node[above] {$c$} (B) -- node[right] {$a$} (C) -- node[below] {$a$} (A);

\draw (1.25cm,-1.0cm) rectangle (1.5cm,-0.75cm);
 \end{tikzpicture}  and so $(a, a, c)$ is a Pythagorean triple.     Now c must either be even or odd.
    \begin{itemize}
    \item
      Suppose $c$ even. Then, by 1., $a$ is even. So by 2., we can double the unit and halve all the dimensions. Clearly, we cannot do this indefinitely, since otherwise the unit will grow larger than $a$.
    \item
      So we must have a Pythagorean triple of the form $(a, a, c)$ in which $c$ is odd. But then, by 4., $a$ is both even and odd, a contradiction. \hfill \mbox{$\Box$}
    \end{itemize}

     The proof of incommensurability is believed to be the first proof by contradiction in the history of mathematical proofs.
  The proof cannot be “seen” by looking at a diagram: it is necessary to follow a sequence of sentences with reasons.

    Theorem~\ref{incommensthe} implies that $\sqrt{2}$ is not a rational number. \\
  Proof: 
    Assume $\sqrt{2} = \myfrac{p}{q}$, then $2q^2 = p^2$.
  Hence $(q, q, p)$ forms a Pythagorean triple. \begin{tikzpicture}[scale=0.25]
      \coordinate  (A) at (-1cm,-1.cm);
\coordinate  (C) at (1.5cm,-1.0cm);
\coordinate  (B) at (1.5cm,1.0cm);
\draw (A) -- node[above] {$p$} (B) -- node[right] {$q$} (C) -- node[below] {$q$} (A);

\draw (1.25cm,-1.0cm) rectangle (1.5cm,-0.75cm);
 \end{tikzpicture}
  Hence there is a unit which measures exactly the side and diagonal of a square.
  This contradicts the incommensurability theorem. \hfill \mbox{$\Box$}
  \subsection{Numbers and Magnitudes}
  \label{numbmagsec}
  With the incommensurability results, the notion of ``number'' as a discrete  collection of units (e.g., naturals or rationals) was no longer enough.  There arose a need for numbers that are  continuous.
   The Greeks did not know how to handle these continuous quantities.
   The main problem was that they treated  mathematical objects as given and did not conceive of constructing them.  And so, they juggled with two notions:
           \begin{itemize}
    \item
      Their notion of ``numbers'' (as a multitude of units, Definition 2 of Book VII).
    \item
      The so-called ``magnitudes'' (which include things like lines and areas and volumes, etc.).
        \end{itemize}
            The Greeks developed arithmetic for their numbers, but treated their magnitudes geometrically.
      However, although they had not thought of constructing new mathematical objects, they did introduce a procedure for approximating ratios. Such approximations were helpful for the much later constructions of magnitudes (e.g., the real numbers).
      
Before explaining how the Greeks developed approximations, we explain the    {\em anthyphairesis} concept.    
Anthyphairesis  is composed of two Greek terms: ${\upsilon\phi\alpha\iota\rho\epsilon\omega}$ (meaning {\em subtract})
and ${\alpha\nu\tau\iota}$ (meaning {\em alternating/reciprocal}) and hence
${\alpha\nu\theta\upsilon\phi\alpha\iota\rho\epsilon\sigma\iota\varsigma}$ stands for {\em alternated/reciprocal subtraction}.
So, given whole numbers $r_0$ and $r_1$, repeatedly subtract $r_1$ from $r_0$, $r_0-r_1$,  $r_0-r_1-r_1$, $\cdots$ until $r_2<r1$ remains, 
then repeat the process for $r_1$ and $r_2$, and so on.

  \begin{tikzpicture}[scale=0.55]
      \path (0,0) coordinate (A) (0,2) coordinate (B) (2,0) coordinate (C) (2,2) coordinate (D);
     \draw (A)
      -- (B) node [midway, left] {$r_1$}
      -- (D) node [midway, above] {$r_1$}
      -- (C)
      --(A)
      -- cycle;
          \path (2,0) coordinate (A1) (2,2) coordinate (B1) (4,0) coordinate (C1) (4,2) coordinate (D1);
     \draw (A1)
      -- (B1) 
      -- (D1) node [midway, above] {$r_1$}
      -- (C1)
      --(A1)
      -- cycle;
             \path (4,0) coordinate (A2) (4,2) coordinate (B2) (6,0) coordinate (C2) (6,2) coordinate (D2);
     \draw (A2)
      -- (B2) 
      -- (D2) node [midway, above] {$r_1$} 
      -- (C2) 
      --(A2)
      -- cycle;
               \path (6,0) coordinate (A3) (6,2) coordinate (B3) (8,0) coordinate (C3) (8,2) coordinate (D3);
     \draw (A3)
      -- (B3) 
      -- (D3) node [midway, above] {$r_1$}
      -- (C3) 
      --(A3)
      -- cycle;
      \path 
      (8.2,0) coordinate (C31) (8.2,2) coordinate (D31);
      \tkzDrawSegments(C3,C31);
      \tkzDrawSegments(D3,D31);
\draw (8.4,1) -- (8.5,1) node {...};
      \path (9,0) coordinate (A4) (9,2) coordinate (B4) (11,0) coordinate (C4) (11,2) coordinate (D4)  (11.5,0) coordinate (C5) (11.5,2) coordinate (D5) (8.8,0) coordinate (C32) (8.8,2) coordinate (D32)(11,0.5) coordinate (Z1)(11.5,0.5) coordinate (Z2)(11,1) coordinate (Z3)(11.5,1) coordinate (Z4)(11,1.5) coordinate (Z5)(11.5,1.5) coordinate (Z6);
\draw(Z1)--(Z2);
\draw(Z3)--(Z4) node[midway, above]{$:$};
\draw(Z5)--(Z6);
     \draw (A4)
      -- (B4) 
      -- (D4) node [midway, above] {$r_1$}
      -- (C4) 
      --(A4)
      -- cycle;
       \tkzDrawSegments(C4,C32);
      \tkzDrawSegments(D4,D32);
      \draw (D4)--(D5) node [midway, above] {$r_2$};
      \draw (C4)--(C5);
       \draw (C5)--(D5);
\draw (4,3) circle [radius=0.4] node {$r_0$};
  \end{tikzpicture}
  
  Euclid used anthyphairesis 
to check whether two numbers are prime to one another.  He proved that anthyphairesis applied to two relatively prime numbers leads to the unit.
\begin{center}   \begin{tcolorbox}[width=.88\textwidth, colframe=red]
    \begin{center}
      \label{propXI.1}
      
    \textsc{Proposition 1.\  of Book VII of the \emph{Elements}}\index{Book VII ! Proposition 1}
  \end{center}
 {\small Two unequal numbers being set out, and the less being continuously subtracted in turn from the greater, if the number left never measures the one before it until a unit is left, the original numbers will be prime to one another.}
\end{tcolorbox}\end{center}

\begin{example}
  Here is why 17 and 3 are prime to one another.
       \begin{center}
        $      \begin{array}{|lllllll|}
        \hline
        17 & = & 5 & \times & 3 & + & 2 \\
        \hline
          3 & = & 1 & \times & 2 & + & \textcircled{1} \\
          \hline
            2 & = & 2 & \times & 1 & + & 0 \\
        \hline
      \end{array}
        $
        \end{center}
       The ratio          and         continued fraction are respectively:\\
       $[5, 1, 2]$ and 
  $\myfrac{17}{3}=\mbox{\fbox{$5$}}
  +\myfrac{1}{\mbox{\fbox{$1$}}+\myfrac{1}{\mbox{\fbox{$2$}}}}.$
            \begin{center}
  \begin{tikzpicture}[scale=0.3]
    \path (-6,0) coordinate (A1) (-3,0) coordinate (A2) (-3,3) coordinate (B2) (-6,3) coordinate (B1) (0,0) coordinate (A) (0,3) coordinate (B) (3,0) coordinate (C)  (6,0) coordinate (N) (6,3) coordinate (M)  (3,3) coordinate (D) (9,0) coordinate (C1) (9,3) coordinate (D1)  (9,1) coordinate (D103)  (11,0) coordinate (C2) (11,3) coordinate (D2)  (10,0) coordinate (C110)
    (11,1) coordinate (C21)  (10,1) coordinate (C22) (9,1) coordinate (C101);
    \draw (A2)
      -- (A1) node [midway, below] {$3$}
      -- (B1) node [midway, left] {$3$}
      -- (B2) 
      -- (A2)
      -- cycle;
      \draw (A)
      -- (A2) node [midway, below] {$3$}
      -- (B2) 
      -- (B)
      -- (A)
      -- cycle;
      \draw (C)
      -- (A) node [midway, below] {$3$}
      -- (B) 
      -- (D)
      -- (C)
            -- cycle;
      \draw (C)
      -- (N) node [midway, below] {$3$}
      -- (C1) node [midway, below] {$3$}
      -- (D1) 
      -- (D);
           \draw (C1)
      -- (C110) node [midway, below] {$1$}
      -- (C2)  node [midway, below] {$1$}
      -- (C21)node [midway, right] {$1$}
      -- (C101);
      \draw (C110)--(C22);
      \draw(C21)--(D2)node [midway, right] {$2$};
      \draw(D2)--(D1)node [midway, above] {$2$};
      \draw (N)--(M);
      \draw (2,3.5) circle [radius=0.6] node {$17$};
  \end{tikzpicture}
  \end{center}
\end{example}
Euclid proved that
anthyphairesis  applied to non relatively prime numbers gives their greatest common divisor (GCD).
\begin{center}   \begin{tcolorbox}[width=.88\textwidth, colframe=red]
  \begin{center}
     \textsc{Proposition 2.\ of Book VII of the \emph{Elements}}\index{Book VII ! Proposition 2}
  \end{center}
 Given two numbers not prime to one another, to find their greatest common measure.
\end{tcolorbox}\end{center}
\begin{example}
  As we see below,     136 and 6 are  not prime to one another and  their greatest common divider       is 2.
  \begin{center}
            $      \begin{array}{|lllllll|}
        \hline
        136 & = & 22 & \times & 6 & + & 4.\\
        \hline
          6 & = & 1 & \times & 4 & + &  \textcircled{2}.\\
          \hline
           4 & = & 2 & \times & 2 & + & 0.\\
        \hline
      \end{array}
    $\end{center}
   The ratio          and         continued fraction are respectively:\\
  $[22, 1, 2]$ and 
  $\myfrac{136}{6}=\mbox{\fbox{$22$}}
   +\myfrac{1}{\mbox{\fbox{$1$}}+\myfrac{1}{\mbox{\fbox{$2$}}}}.$
   \begin{center}
 \begin{tikzpicture}[scale=0.3]
    \path (-6,0) coordinate (A1) (-3,0) coordinate (A2) (-3,3) coordinate (B2) (-6,3) coordinate (B1) (0,0) coordinate (A) (0,3) coordinate (B) (3,0) coordinate (C)  (6,0) coordinate (N) (6,3) coordinate (M)  (3,3) coordinate (D) (9,0) coordinate (C1) (9,3) coordinate (D1)  (9,1) coordinate (D103)  (11,0) coordinate (C2) (11,3) coordinate (D2)  (10,0) coordinate (C110)
    (11,1) coordinate (C21)  (10,1) coordinate (C22) (9,1) coordinate (C101);
    \draw (A2)
      -- (A1) node [midway, above] {$6$}
      -- (B1) node [midway, left] {$6$}
      -- (B2) 
      -- (A2)
      -- cycle;
       \draw (1,2.5)  node {..........................};
       \draw (1,-1) -- (1,-1) node {22 repeats of $6\times6$ square};
      \draw (1,0.5)  node {..........................};
      \draw (1,1.5)  node {..........................};
               \draw (A2) -- (N);
           \draw (B2) -- (M);
            -- cycle;
      \draw (C)
      -- (N) 
      -- (C1) node [midway, above] {$6$}
      -- (D1) 
      -- (D);
           \draw (C1)
      -- (C110) node [midway, below] {$1$}
      -- (C2)  node [midway, below] {$1$}
      -- (C21)node [midway, right] {$1$}
      -- (C101);
      \draw (C110)--(C22);
      \draw(C21)--(D2)node [midway, right] {$4$};
      \draw(D2)--(D1)node [midway, above] {$4$};
      \draw (N)--(M);
      \draw (2,3.5) circle [radius=0.85] node {$136$};
      
 \end{tikzpicture}
 \end{center}
\end{example}
\begin{example}
  \begin{itemize}
    \item
      12 and 5 are prime to one another.
      \begin{center}
     $      \begin{array}{|lllllll|}
        \hline
           12 & = & 2 & \times & 5 & + & 2 \\
            \hline
          5 & = & 2 & \times & 2 & + & \textcircled{1} \\
               \hline
            2 & = & 2 & \times & 1 & + & 0 \\
              \hline
      \end{array}
        $
      \end{center}
       The ratio          and         continued fraction are respectively:\\
   $[2,2,2]$ and 
  $\myfrac{12}{5}=\mbox{\fbox{$2$}}
  +\myfrac{1}{\mbox{\fbox{$2$}}+\myfrac{1}{\mbox{\fbox{$2$}}}}.$
   \begin{center}
  \begin{tikzpicture}[scale=0.25]
      \path (0,0) coordinate (A) (0,5) coordinate (B) (5,0) coordinate (C) (5,5) coordinate (D) (10,0) coordinate (C1) (10,5) coordinate (D1)  (10,3) coordinate (D103)  (12,0) coordinate (C2) (12,5) coordinate (D2)  (11,0) coordinate (C110) (12,1) coordinate (C21) (12,3) coordinate (D21) (11,1) coordinate (C22) (10,1) coordinate (C101);
     \draw (B)
      -- (A) node [midway, left] {$5$}
      -- (C) node [midway, below] {$5$}
      -- (D)
      --(B)
      -- cycle;
         \draw (C)
      -- (C1) node [midway, below] {$5$}
      -- (D1) 
      -- (D);
           \draw (C1)
      -- (C110) node [midway, below] {$1$}
      -- (C2)  node [midway, below] {$1$}
      -- (C21)node [midway, right] {$1$}
      -- (C101);
      \draw (C110)--(C22);
      \draw(C21)--(D21)node [midway, right] {$2$};
      \draw(D21)--(D2)node [midway, right] {$2$};
      \draw(D2)--(D1)node [midway, above] {$2$};
      \draw(D21)--(D103);
\draw (5,5.5) circle [radius=0.7] node {$12$};
      \end{tikzpicture}
   \end{center}
\item
  22 and 6 are  not prime to one another and  their greatest common divider       is 2.
  \begin{center}
            $      \begin{array}{|lllllll|}
        \hline
            22 & = & 3 & \times & 6 & + & 4.\\
            \hline
          6 & = & 1 & \times & 4 & + &  \textcircled{2}.\\
                \hline
            4 & = & 2 & \times & 2 & + & 0.\\
              \hline
      \end{array}
    $
    \end{center}
   The ratio          and         continued fraction are respectively:\\
     $[3, 1, 2]$ and 
       $\myfrac{22}{6}=\mbox{\fbox{$3$}}
   +\myfrac{1}{\mbox{\fbox{$1$}}+\myfrac{1}{\mbox{\fbox{$2$}}}}.$
   \begin{center}
  \begin{tikzpicture}[scale=0.3]
    \path (0,0) coordinate (A) (0,3) coordinate (B) (3,0) coordinate (C)  (6,0) coordinate (N) (6,3) coordinate (M)  (3,3) coordinate (D) (9,0) coordinate (C1) (9,3) coordinate (D1)  (9,1) coordinate (D103)  (11,0) coordinate (C2) (11,3) coordinate (D2)  (10,0) coordinate (C110)
    (11,1) coordinate (C21)  (10,1) coordinate (C22) (9,1) coordinate (C101);
     \draw (B)
      -- (A) node [midway, left] {$6$}
      -- (C) node [midway, below] {$6$}
      -- (D)
      --(B)
      -- cycle;
      \draw (C)
      -- (N) node [midway, below] {$6$}
      -- (C1) node [midway, below] {$6$}
      -- (D1) 
      -- (D);
           \draw (C1)
      -- (C110) node [midway, below] {$2$}
      -- (C2)  node [midway, below] {$2$}
      -- (C21)node [midway, right] {$2$}
      -- (C101);
      \draw (C110)--(C22);
      \draw(C21)--(D2)node [midway, right] {$4$};
      \draw(D2)--(D1)node [midway, above] {$4$};
      \draw (N)--(M);
\draw (5,3.5) circle [radius=0.6] node {$22$};
          \end{tikzpicture}
   \end{center}
    \end{itemize}
\end{example}

The Greeks 
also applied anthyphairesis to magnitudes.    They showed that two magnitudes are commensurable if and only if anthyphairesis terminates and that     if the anthyphairesis procedure of finding the ratio or GCD of two numbers is applied to incommensurable magnitudes, it will not terminate.  
  \begin{center}   \begin{tcolorbox}[width=.88\textwidth, colframe=red]
  \begin{center}
      \textsc{Proposition 2 of Book X of the \emph{Elements}.}
                              \end{center}
   If, when the less of two unequal magnitudes\index{magnitude} is continuously subtracted in turn from the greater, that which is left never measures the one before it, the magnitudes\index{magnitude} will be incommensurable.
  \end{tcolorbox}\end{center}
  \begin{example}
    We show that $\sqrt{2}$ is incommensurable.
    
    \begin{tikzpicture}[scale=0.65]
  \coordinate [label={left:$C$}] (C) at (0, 0);
  \coordinate [label={above:$B$}] (B) at (0, 4);
  \coordinate [label={right:$A$}] (A) at (4, 0);
   \fill (A) circle[radius=2pt];
   \fill (C) circle[radius=2pt];
   \fill (B) circle[radius=2pt];
    \newcommand{\ranglesize}{0.3cm}
  \draw (C) -- ++ (0, \ranglesize) -- ++ (\ranglesize, 0) -- ++ (0, -\ranglesize);
  \tkzDrawPolygon[line width = 1pt](A,B,C)
\tkzLabelSegment[below](C,A){$1$}
\tkzLabelSegment[left](B,C){$1$}
\tkzLabelSegment[above right](B,A){$\sqrt{2}$}
\end{tikzpicture}
\hspace{0.1in}
\begin{tikzpicture}[scale=0.65]
  \coordinate [label={left:$C$}] (C) at (0, 0);
  \coordinate [label={above:$B$}] (B) at (0, 4);
  \coordinate [label={right:$A$}] (A) at (4, 0);
  \coordinate [label={right:$D$}] (D) at (1.35, 2.65);
  \coordinate [label={left:$F$}] (F) at (0, 1.65);
  \fill (D) circle[radius=2pt];
   \fill (A) circle[radius=2pt];
   \fill (C) circle[radius=2pt];
   \fill (B) circle[radius=2pt];
    \fill (F) circle[radius=2pt];
    \draw (D) -- (F);
    \draw [dashed] (F) -- (A);
  \newcommand{\ranglesize}{0.3cm}
  \draw (C) -- ++ (0, \ranglesize) -- ++ (\ranglesize, 0) -- ++ (0, -\ranglesize);
  \tkzMarkRightAngle[draw=blue,size=.2](F,D,A);
\tkzDrawPolygon[line width = 1pt](A,B,C)
\tkzLabelSegment[above right](B,D){$\sqrt{2}-1$}
\tkzLabelSegment[left](C,F){$\sqrt{2}-1$}
\tkzLabelSegment[above right](D,A){$1$}
\tkzLabelSegment[left](B,F){$2-\sqrt{2}$}
\end{tikzpicture}

  \begin{tikzpicture}[scale=0.65]
  \coordinate [label={left:$C$}] (C) at (0, 0);
  \coordinate [label={above:$B$}] (B) at (0, 4);
  \coordinate [label={right:$A$}] (A) at (4, 0);
  \coordinate [label={right:$D$}] (D) at (1.35, 2.65);
  \coordinate [label={right:$G$}] (G) at (0.8, 3.2);
  \coordinate [label={left:$F$}] (F) at (0, 1.65);
  \coordinate [label={left:$E$}] (E) at (0, 3.2);
  \fill (D) circle[radius=2pt];
  \fill (G) circle[radius=2pt];
  \fill (E) circle[radius=2pt];
   \fill (A) circle[radius=2pt];
   \fill (C) circle[radius=2pt];
   \fill (B) circle[radius=2pt];
    \fill (F) circle[radius=2pt];
    \draw (D) -- (F);
     \draw (E) -- (G);
     \draw [dashed] (F) -- (A);
      \draw [dashed] (F) -- (G);
  \newcommand{\ranglesize}{0.3cm}
  \draw (C) -- ++ (0, \ranglesize) -- ++ (\ranglesize, 0) -- ++ (0, -\ranglesize);
  \tkzMarkRightAngle[draw=blue,size=.2](F,D,A);
  \tkzMarkRightAngle[draw=blue,size=.2](F,E,G);
\tkzDefSquare(B,A) \tkzGetPoints{E}{F}
\tkzDefSquare(A,C) \tkzGetPoints{G}{H}
\tkzDefSquare(C,B) \tkzGetPoints{I}{J}
\tkzDrawPolygon[line width = 1pt](A,B,C)
\coordinate [label={left:$3-2\sqrt{2}$}] (M) at (-0.25, 3.75);
\coordinate [label={left:$\sqrt{2}-1$}] (N) at (-0.25, 1);
\coordinate [label={left:$\sqrt{2}-1$}] (O) at (-0.25, 2.75);
  \end{tikzpicture}
  \hspace{0.1in}
\begin{tikzpicture}[scale=0.65]
  \coordinate [label={left:$C$}] (C) at (0, 0);
  \coordinate [label={above:$B$}] (B) at (0, 4);
  \coordinate [label={right:$A$}] (A) at (4, 0);
  \coordinate [label={right:$D$}] (D) at (1.35, 2.65);
  \coordinate [label={right:$G$}] (G) at (0.8, 3.2);
   \coordinate [label={right:$H$}] (H) at (0.3, 3.7);
  \coordinate [label={left:$F$}] (F) at (0, 1.65);
  \coordinate [label={left:$E$}] (E) at (0, 3.2);
   \coordinate [label={left:$I$}] (I) at (0, 3.6);
   \fill (D) circle[radius=2pt];
   \fill (H) circle[radius=2pt];
   \fill (I) circle[radius=2pt];
  \fill (G) circle[radius=2pt];
  \fill (E) circle[radius=2pt];
   \fill (A) circle[radius=2pt];
   \fill (C) circle[radius=2pt];
   \fill (B) circle[radius=2pt];
    \fill (F) circle[radius=2pt];
    \draw (D) -- (F);
     \draw (I) -- (H);
     \draw (E) -- (G);
     \draw [dashed] (F) -- (A);
     \draw [dashed] (G) -- (I);
      \draw [dashed] (F) -- (G);
  \newcommand{\ranglesize}{0.3cm}
  \draw (C) -- ++ (0, \ranglesize) -- ++ (\ranglesize, 0) -- ++ (0, -\ranglesize);
  \tkzMarkRightAngle[draw=blue,size=.2](F,D,A);
  \tkzMarkRightAngle[draw=blue,size=.2](F,E,G);
  \tkzMarkRightAngle[draw=blue,size=.1](B,H,I);
\tkzDefSquare(B,A) \tkzGetPoints{E}{F}
\tkzDefSquare(A,C) \tkzGetPoints{G}{H}
\tkzDefSquare(C,B) \tkzGetPoints{I}{J}
\tkzDrawPolygon[line width = 1pt](A,B,C)
\end{tikzpicture}

    Geometrically, we assume the isosceles rectangular triangle $BCA$ below and take $BD$ of length $\sqrt{2}-1$.  From $D$ we draw the perpendicular to $AB$ meeting $BC$ on $F$. We get an isosceles recltangular triangle $BDF$.  We repeat the process obtaining  isosceles recltangular triangles $BEG$, $BHI$, and so on.  In this repetition, the less of two unequal magnitudes is continuously subtracted in turn from the greater, yet what is left never
measures the one before it.  This  can be repeated infinitely and $\sqrt{2}$ is incommensurable. Using anthyphairesis: 
  \begin{center}
    \begin{tikzpicture}[scale=1.5]
    \path (0,0) coordinate (C) (0,1) coordinate (D) (1,0) coordinate (C1) (1,1) coordinate (D1)  (1,.586) coordinate (D103)  (1.414,0) coordinate (C2) (1.414,1) coordinate (D2)  (1.242,0) coordinate (C110) (1.07,0) coordinate (C1101)
    (1.414,.172) coordinate (C21) (1.414,.586) coordinate (D21)
    (1.242,.172) coordinate (C22)(1.07,.172) coordinate (D1101)
    (1,.172) coordinate (C101);
    \draw(C22)--(C110);
    \draw(D1101)--(C1101) node [midway, left] {$\ddots$};
    \draw(D1101)--(C1101) node [midway, left] {$\vdots$};
         \draw (C)
           -- (D)node [midway, left] {$1$};
              \draw (C)
      -- (C1) node [midway, below] {$1$}
      -- (D1) 
      -- (D);
           \draw (C1)
      -- (C110) 
      -- (C2)  
      -- (C21)node [midway, right] {$3-2\sqrt{2}$};
      -- (C101);
      \draw (C101)--(C21);
      \draw(C21)--(D21)node [midway, right] {$\sqrt{2}-1$};
      \draw(D21)--(D2)node [midway, right] {$\sqrt{2}-1$};
      \draw(D2)--(D1)node [midway, above] {$\sqrt{2}-1$};
      \draw(D21)--(D103);
\draw (0.5,1.2) circle [radius=0.2] node {$\sqrt{2}$};
  \end{tikzpicture}
  \end{center}
  The ratio of $\sqrt{2}$ to 1 is  $[1,2,2,\cdots]$ and   
    the continued fraction is $\sqrt{2} =  1 +\myfrac{1}{2+\myfrac{1}{2+\ddots}}.$

    $\sqrt{2}$ is called a quadratic irrational because it is the solution to the quadratic equation $x^2-2 = 0$.  Note that these continued fractions  provide an approximation to $\sqrt{2}$ as follows:
  \vspace{-0.1in}
\begin{itemize}
  \item
    $\sqrt{2} \approx 1$,
    \item
      $\sqrt{2} \approx 1+\myfrac{1}{2} = 1.5$,
      \item
        $\sqrt{2} \approx 1+\myfrac{1}{2+\myfrac{1}{2}} = 1.4$,
        \item
          $\sqrt{2} \approx 1+\myfrac{1}{2+\myfrac{1}{2+\myfrac{1}{2}}} = 1.417$,
          \item
            $\sqrt{2} \approx 1+\myfrac{1}{2+\myfrac{1}{2+\myfrac{1}{2+\myfrac{1}{2}}}} = 1.4139$ etc.
             \end{itemize}
\end{example}
  Infinite repetitions/approximations were a useful part of Greek's Mathematics but,  anthyphairesis had its limitations.  E.g., the obvious theorem below cannot be proved with it:
  
   \begin{center}   \begin{tcolorbox}[width=.88\textwidth, colframe=red]
    {\it If the ratio of $A$ to $C$ is the same as the ratio of $B$ to $C$, then $A = B$}.
   \end{tcolorbox}\end{center}
   To overcome the problems, Eudoxus, defined {\it proportion} (having the same ratio) for magnitudes instead of ratios.  He invented the method of exhaustion which was used by Archimedes and Euclid (see Sections~\ref{euclidsquaresec} and~\ref{archimesec}). Theodorus of Cyrene used Eudoxus approximation  in his spiral of irrational numbers  pictured earlier.
   \subsection{The Greeks' problems with infinitesimals/limits}
   The Greeks were puzzled by limits and  infinitesimals.  They needed approximations but faced obstacles they could not explain.  For example, in the diagram below,  the length of the stepped line is clearly $2s$ no matter how many steps there
are. But as the number of steps increases, the stepped line seems to approach the
diagonal whose length  is
$\sqrt{2} s \not = 2s$.
\begin{center}
\begin{tikzpicture}[scale=0.5]
\path (0,0) coordinate (A) (0,6) coordinate (B) (6,0) coordinate (C)  (6,6) coordinate (D);
\draw (A)
--(B)
--(D)
--(C)
--(A) node [midway, below] {$s$}
--cycle;
\draw (3,3)--(6,3);
\draw (3,0)--(3,3);

\draw (2.25,1.5)--(2.25,2.25);
\draw (2.25,2.25)--(3,2.25);
\draw (4.5,3)--(4.5,4.5);
\draw (4.5,4.5)--(6,4.5);
\draw (.375,0)--(.375,.375);
\draw (.375,.375)--(.75,.375);

\draw (.375,0)--(.375,.375);
\draw (.375,.375)--(.75,.375);
\draw (1.125,.75)--(1.125,1.125);
\draw (1.125,1.125)--(1.5,1.125);

\draw (1.875,1.5)--(1.875,1.875);
\draw (1.875,1.875)--(2.25,1.875);
\draw (2.625,2.25)--(2.625,2.625);
\draw (2.625,2.625)--(3,2.625);

\draw (3.375,3)--(3.375,3.375);
\draw (3.375,3.375)--(3.75,3.375);
\draw (4.125,3.75)--(4.125,4.125);
\draw (4.125,4.125)--(4.5,4.125);

\draw (4.875,4.5)--(4.875,4.875);
\draw (4.875,4.875)--(5.25,4.875);
\draw (5.625,5.25)--(5.625,5.625);
\draw (5.625,5.625)--(6,5.625);

\draw (3.75,3)--(3.75,3.75);
\draw (3.75,3.75)--(4.5,3.75);

\draw (5.25,4.5)--(5.25,5.25);
\draw (5.25,5.25)--(6,5.25);

\draw (1.5,0)--(1.5,1.5);
\draw (1.5,1.5)--(3,1.5);
\draw (0.75,0)--(0.75,0.75);
\draw (0.75,0.75)--(1.5,0.75);
\draw (A)--(D)[thick];
\end{tikzpicture}
\end{center}

They demonstrated many paradoxes like the following:

 \begin{center}   \begin{tcolorbox}[width=.88\textwidth, colframe=orange]
     {\bf  Zeno's Dichotomy Paradox } There is no motion, because what moves must arrive at the middle of its course before it reaches the end.
 \end{tcolorbox}\end{center}
 
  For example, to leave the room, you first have to get halfway to the door, then halfway from that point to the door, etc.  No matter how close you are to the door, you have to go half the remaining distance. Hence, there is no finite motion because always going half way while in motion is infinite. 

     Suppose the distance is 1 meter and the object  moves at 1 meter per second.
  It must reach  halfway ($\frac{1}{2}$ meter from the starting point) in $a_1 = \frac{1}{2}$ second.  Let  $t_1 = a_1$.
  From this halfway point, the object moves halfway to the end, which is $a_2=\frac{1}{4}$ meters.  The total time  so far is  $t_2 = a_1+a_2= \frac{1}{2} + \frac{1}{4}$.
We clearly have the following infinite sequences:
$$a_1, a_2, a_3, \ldots = \frac{1}{2}, \frac{1}{4}, \frac{1}{8}, \ldots$$
$$t_1, t_2, \ldots = \frac{1}{2}, \frac{3}{4}, \frac{7}{8}, \ldots \mbox{ where each } t_n = a_1 + a_2 +\cdots + a_n.$$
  Zeno concluded that the total time which is the sum of an infinite sequence must be infinite and  we can never reach our destination.
     This is incorrect since we can reach our destination in a finite time.  So, where did Zeno get it wrong?

In modern notation, we see that:
\begin{itemize}
  \item
$t_n = \frac{2^n-1}{2^n} =
  1 - \frac{1}{2^n} <1 \mbox{      and      } \lim _{n \to \infty} t_n = 1 .$
  \item
  $2\Sigma_{n=1}^\infty a_n=  2a_1+ 2\Sigma_{n=2}^\infty \frac{1}{2^n}=
  1+ \Sigma_{n=1}^\infty \frac{1}{2^{n}} =1+ \Sigma_{n=1}^\infty a_n.$
\item
 $ \mbox{Hence, } \Sigma_{n=1}^\infty a_n= 1 \mbox{ and } 
lim_{n\to \infty} t_n = \Sigma_{n=1}^\infty a_n= 1.$
\end{itemize}

Despite the complications of limits, the Greeks continued to use them to measure magnitudes.  Both Archimedes and Euclid (see Sections~\ref{euclidsquaresec} and~\ref{archimesec}) used  Eudoxus theory of proportions which is a geometric method based on exhaustive approximations designed to overcome  the difficulties obtained from the discovery
of the irrationals.

      \subsection{The area  a regular polygon}
      For both Archimedes’ theorem and Euclid’s theorem, we need a general formula for the area of a regular polygon (i.e., a polygon where all angles (resp.\ all sides) are equal). Let us start with the area of a square of side $s$.
\begin{center}
\begin{tikzpicture}[scale=0.25]
\path (0,0) coordinate (A) (0,6) coordinate (B) (6,0) coordinate (C)  (6,6) coordinate (D);
\draw (A)
--(B)
--(D)
--(C)
--(A) node [midway, below] {$s$}
--cycle;
\draw (0,3)--(6,3);
\draw (3,0)--(3,6)  node [pos =.25, right]{h};
\draw (B)--(C)[thick];
\draw (A)--(D)[thick];
\end{tikzpicture}
\end{center}
Instead of simply taking $s^2$, take the bottom of the 4 triangles obtained by the diagonals.  Note that the altitude $h = \myfrac{1}{2}s$.   The area $A$ of the square = 4$\times$ area of triangle = $4\times \frac{1}{2}hs = \frac{1}{2}h(4s) = \frac{1}{2}hp$. where $p$ is the perimeter of the square.\\
Note that  $A = \frac{1}{2}hp =  \frac{1}{2}\frac{s}{2}(4s)= s^2$.

   \begin{center}
     \includegraphics[scale = 0.45]{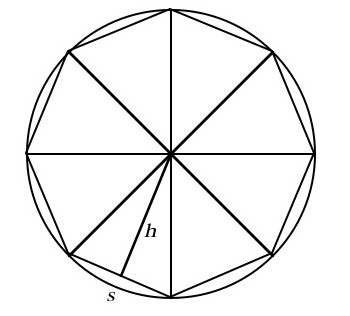}
       \end{center}

Now let us consider a regular octagon. If we divide it into triangles the same way,
we get eight triangles, each of whose areas is $\frac{1}{2}hs$.
 If we take all eight triangles and
note that here $p = 8s$, we get for the area
$A =
\frac{1}{2}h(8s) =\frac{1}{2}hp$.

We saw this for the square and the regular octagon, but it holds for every regular polygon:
\begin{center}   \begin{tcolorbox}[width=.88\textwidth, colframe=red]
    The area of any regular polygon is one-half the altitude to a side times the perimeter, or $\frac{1}{2}hp$.
   \end{tcolorbox}\end{center} 

Now we come to the area of a circle.  Note that the above polygon was inscribed in  the circle with circumference $C$.
     If we  keep increasing the number of sides, the perimeter will approach the circumference $C$ and the altitude will approach the radius $r$.  By the above, this suggests that the formula for the area of a circle should be
\[
A = \frac{1}{2}rC .
\]%
  And since $\pi$ is defined to be the ratio of the circumference of a circle to twice its radius, we have
\[
\pi = \frac{C}{2r} ,
\]%
   Hence 
\[
A = \frac{1}{2}r (2 \pi r) = \pi r^2 
\]%
      This must have seemed obvious to the ancient Greeks from an early
    period in the history of their geometry.   But how could they prove it?    At one time some of them argued that a circle is a regular polygon with infinitely many sides, but they eventually decided that this kind of reasoning is not immune to attacks by sophists.   For just because regular polygons with an increasing number of sides seems to be approaching a circle, does not automatically justify in deducing this formula for the area of a circle.   They found evidence like this to be misleading.  Recall the stepped line which wrongly gave the impression that $\sqrt{2} s = 2s$.

    \subsection{Euclid on Areas of Circles and Squares}
    \label{euclidsquaresec}
      It took a long time for the proof that $A= \frac{1}{2}rC$ to be given.  Although this was obvious to the Greeks, a proof was hard to find.
      Before that proof was given (by Archimedes), Euclid  proved that the areas of circles have the same proportion  as the squares on their diameters (Proposition 2 of Book XII of \emph{Elements}).
      The proof uses Proposition 1 of Book XII.
\begin{center}   \begin{tcolorbox}[width=.88\textwidth, colframe=red]
    \begin{center}
      \label{propXII.1}
    \textsc{Proposition 1  of Book XII of the \emph{Elements}.}\index{Book XII ! Proposition 1}
  \end{center}
  Similar polygons inscribed in circles are to one another as the squares on the diameters of the circles.
\end{tcolorbox}\end{center}

  Similar figures are those which have the same shape.  In similar polygons the corresponding angles are equal and the corresponding sides all have the same proportion.
   
 \begin{center}   \begin{tcolorbox}[width=.88\textwidth, colframe=red]
  The areas $A$ of similar polygons are proportional to:
    \begin{itemize}
      \item
   The squares of their altitudes $h$. 
\item  
  The squares of their perimeters $p$.
\item
  The squares of any of their linear parts.
    \end{itemize}
    $\myfrac{p_1}{p_2} = \myfrac{h_1}{h_2}$ and $\myfrac{A_1}{A_2} =
    \myfrac{h_1}{h_2} \myfrac{h_1}{h_2} = \myfrac{h_1^2}{h_2^2}= \myfrac{p_1^2}{p_2^2}$.
  \end{tcolorbox}\end{center}
 The proof of Proposition 1 of Book XII uses the above and the fact that $AGB$ is similar to $A'G'B'$ below and hence $(\myfrac{AB}{A'B'})^2=(\myfrac{AG}{A'G'})^2=\myfrac{A_1}{A_2} $.
 \begin{center}
  \begin{tikzpicture}[thick, scale=0.4]
        \coordinate (O) at (2,2);
       \def\radius{3cm}
       \draw (O) circle[radius=\radius];
  \fill (O) circle[radius=2pt] node[below left] {O};

      \tkzDefPoint[label=right:{$B$}](5,2){B};
      \tkzDefPoint[label=left:{$A$}](-1,2){A};
      \tkzDefPoint[label=above:{$D$}](2,5){D};
      \tkzDefPoint[label=left:{$C$}](0.05,4.25){C};
      \tkzDefPoint[label=right:{$E$}](3.95,4.25){E};
            \tkzDefPoint[label=right:{$F$}](3.5,-0.55){F};
      \tkzDefPoint[label=below:{$G$}](0.4,-0.5){G};
      \draw (A)--(B);
      \draw [color = red] (G)--(B);
      \draw [color = red] (A)--(F);
      \draw (D)--(E);
      \draw (E)--(F);
       \draw (G)--(F);
      \draw (A)--(C);
      \draw (D)--(C);
       \draw (A)--(G);
 \end{tikzpicture}
 \begin{tikzpicture}[thick, scale=0.3]
        \coordinate (O) at (2,2);
       \def\radius{3cm}
       \draw (O) circle[radius=\radius];
  \fill (O) circle[radius=2pt] node[below left] {$O'$};

      \tkzDefPoint[label=right:{$B'$}](5,2){B};
      \tkzDefPoint[label=left:{$A'$}](-1,2){A};
      \tkzDefPoint[label=above:{$D'$}](2,5){D};
      \tkzDefPoint[label=left:{$C'$}](0.05,4.25){C};
      \tkzDefPoint[label=right:{$E'$}](3.95,4.25){E};
       \tkzDefPoint[label=right:{$F'$}](3.5,-0.55){F};
      \tkzDefPoint[label=below:{$G'$}](0.4,-0.5){G};
      \draw (A)--(B);
      \draw [color = red] (G)--(B);
        \draw [color = red] (A)--(F);
      \draw (D)--(E);
      \draw (E)--(F);
       \draw (G)--(F);
      \draw (A)--(C);
      \draw (D)--(C);
       \draw (A)--(G);
  \end{tikzpicture}
\end{center}
 Now we look at Euclid's proposition 2 and its proof:
  \begin{center}   \begin{tcolorbox}[width=.88\textwidth, colframe=red]
    \begin{center}
      \label{propXII.2}
      \noindent
       \textsc{Proposition 2  of Book XII of the \emph{Elements}.}\index{Book XII ! Proposition 2}
    \end{center}
     Circles are to one another as the squares on the diameters.
\end{tcolorbox}\end{center}

Euclid starts his proof as follows:

\begin{center}   \begin{tcolorbox}[width=.88\textwidth, colframe=red]
      \noindent
   Let $ABCD$, $EFGH$ be circles, and $BD$, $FH$ their diameters; I say that, as the circle $ABCD$ is to the circle $EFGH$, so is the square on $BD$ to the square on $FH$.

  \medskip
  
   \begin{center} 
   \includegraphics[width=2.7in]{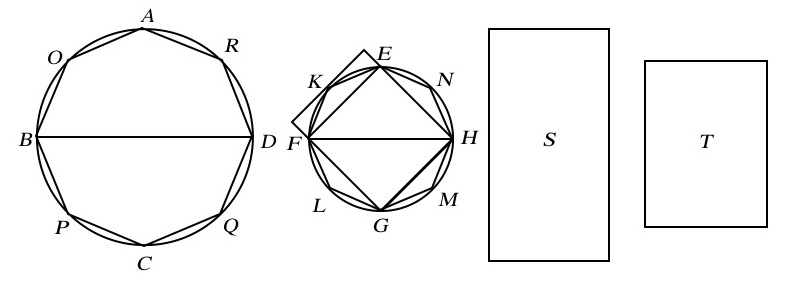} 
   \end{center}
{\small   
For, if the square on $BD$ is not to the square on $FH$ as the circle $ABCD$ is to the circle $EFGH$,
then, as the square on $BD$ is to the square on $FH$, so will the circle $ABCD$ be either to some less area than the circle $EFGH$ or to a greater.}
\end{tcolorbox}\end{center}

Euclid's strategy is to prove his result by contradiction.  In fact, it will be a double proof by contradiction.  He will first assume that it will be in the ratio to a smaller area $S$, derive a contradiction from that, then assume that it will be in the ratio to a larger area $S$, and then derive a contradiction from that area as well.  As a result, the only possibility left will be the result stated in the proposition.

We will not repeat the proof here (see~\cite{Heath:BEE-1956,Kamareddine-seldin-analysis}).  We must mention however that Euclid's method is based on Eudoxus exhaustion which infinitely inscribes and circumscribes polygons inside the circles.   First, Euclid assumes  it to be in that ratio to a less area $S$ and 
 shows that the square $EFGH$ {\em inscribed} in the circle $EFGH$ is greater than half of the circle $EFGH$.   He shows this by noting that the {\em circumscribed} square, which includes area outside the circle, has twice the area of the inscribed square.

 Then, he bisects the circumference $EF$, $FG$, $GH$, $HE$ at the points $K$, $L$, $M$, $N$  and joins $EK$, $KF$, $FL$, $LG$, $GM$, $MH$, $HN$, $NE$
 and proves that the new circumference  (in effect inscribing a new regular polygon with twice the number of sides as the previous one), is  more than half the area inside the circle but outside the previous polygon.  By bisecting the remaining circumferences and joining straight lines, and by doing this continually, one is left with some segments of the circle which will be less than the excess by which the circle EFGH exceeds the area S.

 In modern notation, 
  let the circles have areas $a$ and $b$ respectively, and let the ratio of the squares of their diameters be $k$.
    Let the areas of the polygons inscribed in the circle with area $a$ (resp. $b$) have areas $a_1, a_2, \ldots$ (resp.\  $b_1, b_2, \ldots$).    We have
$
0 < a_1 < a_2 < \ldots < a_n < \ldots < a$
\mbox{ and }
$0 < b_1 < b_2 < \ldots < b_n < \ldots < b$.
\begin{itemize}
\item
For each $n$, we have
\begin{itemize}
\item
$k = \myfrac{a_n}{b_n},
\mbox{ so that }
\myfrac{a_n}{k} = b_n .$
\item
$
  (a - a_{n+1}) < \myfrac{1}{2}(a - a_n)$ and $(b - b_{n+1}) < \myfrac{1}{2}(b - b_n).$
  \end{itemize}
\item
  We want to prove
$k = \myfrac{a}{b} .$
    \item
If $k \not= \myfrac{a}{b}$, then $k = \myfrac{a}{S}$, where $S < b$ or $S > b$.
\begin{itemize}
    \item  Suppose $S < b$.  Choose $N$ so that $b - b_N < b - S .$
The number $N$ represents the number of times the number of sides of the inscribed polygon was doubled.  Then
$S < b_N .$  But
$S = \myfrac{a}{k} > \myfrac{a_N}{k} = b_N ,$
a contradiction.
\item
 Suppose $S > b$.  This is similar to the above case  with $a$ and $b$ reversed.
\end{itemize}
It follows that
$k = \myfrac{a}{b} . \hfill \mbox{$\Box$}
$
\end{itemize}

\subsection{Archimedes' Measurement of a Circle}
\label{archimesec}
  Archimedes used Eudoxus' exhaustion to prove the following proposition (and hence its corollary that the area of a circle of circumference $C$ and radius $r$ is     $A= \myfrac{1}{2}rC$). 
  
  \begin{center}   \begin{tcolorbox}[width=.88\textwidth, colframe=red]
  \begin{center}
     \textsc{Proposition 1  of Archimedes's Book ``Measurement of a Circle''.}
\end{center}

\emph{The area of any circle is equal to a right-angled triangle in which one of the sides about the right triangle is equal to the radius, and the other to the circumference of the circle.}
    \end{tcolorbox}
  \end{center}
As we see from the begin of its proof, an infinite number of polygons will be inscribed/circumscribed in the circle.
   \begin{center}   \begin{tcolorbox}[width=.88\textwidth, colframe=red]
Let $ABCD$ be the given circle, $K$ the triangle described.
\begin{center}
  \includegraphics[width=1.5in]{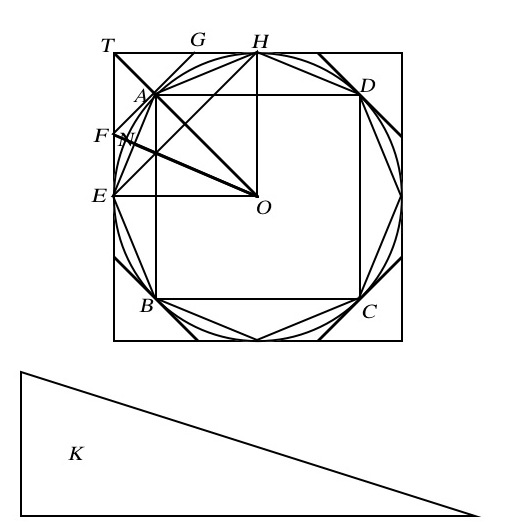} 
\end{center}
Then, if the circle is not equal to $K$, it must be either greater or less.

I.  If possible, let the circle be greater than $K$.

Inscribe a square $ABCD$, bisect the arcs $AB$, $BC$, $CD$, $DA$, and then bisect (if necessary) the halves, and so on, until the sides of the inscribed polygon whose angular points are the points of division subtend segments whose sum is less than the excess of the circle over $K$.
\end{tcolorbox}\end{center}
 
   Let us write the proof in modern notation.
  
  Let $K = \myfrac{1}{2}rC$ (the area of the triangle).  If $A \not= K$, then:
\begin{itemize}
  \item[I.]  Suppose  $A > K$. 
    \begin{itemize}
    \item
        Inscribe a square with side $s_1$, altitude to the side  $h_1$, and  perimeter $p_1$.  The area of the square is $a_1 = \myfrac{1}{2}h_1p_1 .$
  \item
    Now, double the number of sides of the inscribed polygon, and keep doubling it.  For polygon $n$ with side  $s_n$,  altitude to the side  $h_n$, and perimeter   $p_n$,  the area is $a_n = \myfrac{1}{2}h_np_n .$
    \item
From the geometry of the situation, we have that\\
$
h_1 < h_2 < \ldots < h_n < \ldots r,
$\\
$
p_1 < p_2 < \ldots < p_n < \ldots < C,
$
and\\
$
a_1 < a_2 < \ldots < a_n < \ldots < A.
$
\item
Now choose $N$ so that
$
A - a_N < A - \myfrac{1}{2}rC.
$
It follows that
$
\myfrac{1}{2}rC < a_N .
$
\item
But since $h_N < r$, $p_N < C$, and $a_N = \myfrac{1}{2}h_Np_N$, we have
$
a_N < \myfrac{1}{2}rC ,
$
a contradiction.
\end{itemize}
  \item[II.]
  Suppose, on the contrary, that $A < K$.
   \begin{itemize}
   \item
     Circumscribe a square with perimeter $P_1$; then the area is
$A_1 = \myfrac{1}{2}rP_1 .$
\item
Double the number of sides of the circumscribed figure, and keep doing it.  If, for the $n$th polygon, the perimeter is $P_n$, then the area is
$
A_n = \myfrac{1}{2}rP_n .
$
\item
From the geometry, we have\\
$
C < \ldots < P_n < \ldots < P_2 < P_1$
and\\
$
A < \ldots A_n < \ldots < A_2 < A_1 .
$
\item
Choose $N$ where
$
A_N - A < \myfrac{1}{2}rC - A .
$\\
Then
$
A_N < \myfrac{1}{2}rC .
$
\item
But $C < P_N$ and $A_N = \myfrac{1}{2}rP_N$, so
$
\myfrac{1}{2}rC < A_N ,
$
another contradiction.
   \end{itemize}
   \end{itemize}
  It follows that $A = K = \myfrac{1}{2}rC$. \hfill \mbox{$\Box$}
\subsection{Eudoxus, the infinitesimal and the limit}
  We saw the use of Eudoxus' exhaustion method in the proofs of Euclid and Archimedes.  This method infinitely constructs new objects that would eventually only differ in infinitesimal amounts.   It can  be used to develop a  definition  of the limit of a sequence and a function.       Historically, the development of calculus and analysis in European mathematics occured before a definition of the real numbers.  At the time of Descartes, Leibniz and Newton, it had not even been settled whether or not there were infinitely small quantities.   For centuries before and after, infinitesimals oscillated between being accepted and being rejected.  
  They were
                     \emph{introduced} in 450 BC,
                  \emph{banned} by Eucledian mathematicians because 
          of  the problems they faced with them, used by Kepler to calculate the area of an ellipse as the infinite sum of vertical lines contained in the ellipse, 
                  \emph{banned again}
          in  the 1630s
          by religious clerics in Rome.
                        They still \emph{flourished}
                        in the 17th century\footnote{In the ideas that a curved line is made of infinitely small straight line segments, and
                    quantities that differ by an infinitely small quantity are  equal.} and 
                            were \emph{crucial} for the development of calculus by Newton and Leibniz. They were {\bf thought} to exist by Cauchy who used them in his approach to calculus, then they were
                      \emph{abandoned again} in the 19th century due to their unclear logical status to be                   \emph{revived again} in the 20th century especially in Robinson's non-standard analysis.  Nowadays, they take center stage in the foundations of mathematics which many people define as a sound theory of infinitesimals.

The next graph demonstrates how a curved line is made of infinitely small straight line segments.
\begin{center}
  \begin{tikzpicture}[scale = 0.65]
     \tkzInit[xmax=6,ymax=6,xmin=0,ymin=0]
     \draw[red] (0,6) -- (1,0);
     \draw[red] (0,5) -- (2,0);
     \draw[red] (0,4) -- (3,0);
     \draw[red] (0,3) -- (4,0);
     \draw[red] (0,2) -- (5,0);
     \draw[red] (0,1) -- (6,0);
  \end{tikzpicture}
  \hspace{0.3in}
\begin{tikzpicture}[scale = 0.65]
     \tkzInit[xmax=6,ymax=6,xmin=0,ymin=0]
     \draw[red] (0,6) -- (0.5,0);
     \draw[red] (0,5.5) -- (1,0);
     \draw[red] (0,5) -- (1.5,0);
     \draw[red] (0,4.5) -- (2,0);
    \draw[red] (0,4) -- (2.5,0);
     \draw[red] (0,3.5) -- (3,0);
     \draw[red] (0,3) -- (3.5,0);
     \draw[red] (0,2.5) -- (4,0);
     \draw[red] (0,2) -- (4.5,0);
     \draw[red] (0,1.5) -- (5,0);
     \draw[red] (0,1) -- (5.5,0);
      \draw[red] (0,0.5) -- (6,0);
  \end{tikzpicture}
\end{center}
The next example explains a cleric position on infinitesimals.
\begin{example}
                To find the derivative  $f'(2)$ at $x = 2$ of $y = f(x) = x^2$, we assume  $x \not = 2$.  Then we calculate:
          \[
         \myfrac{\Delta y}{\Delta x} = \myfrac{f(x)-f(2)}{x-2} = \myfrac{x^2-2^2}{x-2} =\myfrac{(x+2)(x-2)}{x-2} = x+2.\]
  Since we are only able to conclude that the quotient is equal to x + 2 on the assumption that $x \not= 2$, we appear to have taken an illegal step.
        We justify this by saying that we are
           taking its limit as $x \rightarrow 2$ and write:
      $\myfrac{dy}{dx} = \lim_{x\to 2}\myfrac{\Delta y}{\Delta x}.$
         Newton calls 
    $\myfrac{dy}{dx} = \lim_{x\to 2}\myfrac{\Delta y}{\Delta x},$
         \emph{ultimate value} or \emph{value at instant of disappearance}.  Sarcastically, this is called  \emph{the ghosts of
          a departed quantity} 
in a critique~\cite{berkeley:Anls} by Bishop Berkeley addressed to a certain ``Infidel Mathematician''. \cite{berkeley:Anls}  examined whether the object and principles of the modern Analysis are more distinctly conceived, or more evidently deduced, than religious mysteries and points of faith. 
                             \end{example}
\subsection{Infinitesimals and the birth of analysis}
 At school,  after studying arithmetic and elementary algebra, you are introduced to geometry (\emph{the study of shapes})  and trigonometry (\emph{the study of side lengths and angles of triangles}) and then  
              you move to a  \emph{pre-calculus course} which  combines advanced algebra and geometry with trigonometry.
                  After all this, you are introduced to  \emph{calculus}.
                  Calculus (originally called \emph{infinitesimal calculus}) is the mathematical study of continuous change. The infinitesimal part is important.  
                  It is believed that if Descartes had expressed rather than supressed the infinitesimals and infinites  in his method, he would have invented the calculus before Newton and Leibniz.

                  Calculus formalizes the study of continuous change, while analysis provides it with a rigorous foundation in logic.  As we saw, the Greeks dealt with discrete numbers arithmetically and with continuous magnitudes geometrically.  But continuous systems can be subdivided indefinitely, and their description requires the real numbers.
                  This infinite subdivision  was influenced by  Eudoxus' and Archimedes' approximations. The real numbers   were not present in   the historic approach to define limits and develop the calculus.              The ancient Greeks separated whole and rational numbers, which are discrete, from continuous magnitudes.  They  had different kinds of magnitudes for lengths, areas, volumes, angles, etc., and  never multiplied two lengths to get another length.  The beginning of algebra and 
                  the reduction of  geometrical problems into algebraic and arithmetical ones in  the 9th century~\cite{rashed:arabmaths,rashed:arabmathsfrench} paved the way for  Descartes   innovative ruler-and-compass construction for multiplying two lengths to get a length.
       This  allowed Algebra to be a science concerned with numbers rather than geometric magnitudes.
  
       Here is how the ruler-and-compass construction works:
       \begin{example}
  The length of $AB$ is $a$.  On a line $AC$ through $A$ and at an angle to $AB$, let the length of $AC$ be a unit, and construct $E$ on the same line so that the length of $AE$ is $b$.  Join $C$ and $B$ with line segment $BC$, and construct a line through $E$ parallel to $BC$; let this line intersect the extension of $AB$ at $D$.  Then triangles $ABC$ and $ADE$ are similar. Hence, 
$AE$ is to $AC$ as $AD$ is to $AB$. I.e., 
      $AD = ab$.
  
       \begin{center}
                  \includegraphics[width=1.9in]{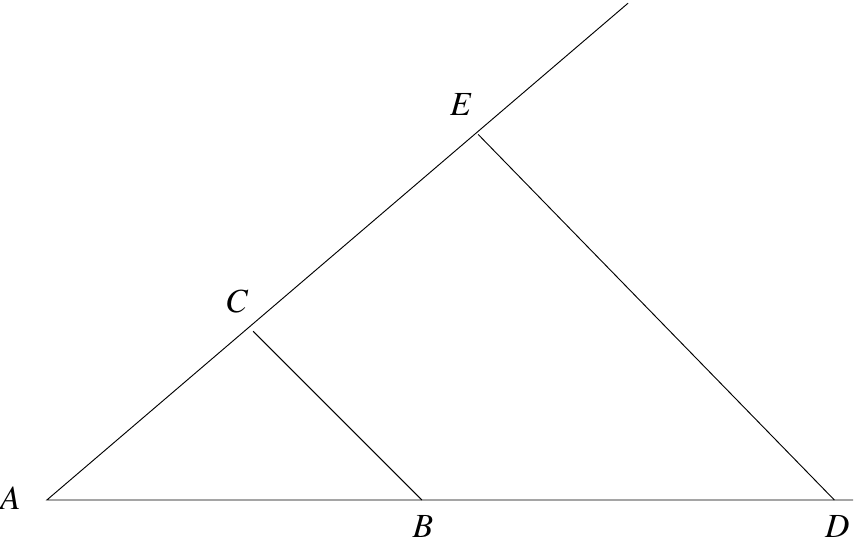}
         \end{center}
\end{example}

        The  move to  generalise the geometric concepts and methods of the calculus to more algebraic forms continued into the 18th century. 
       But the field was still rife with disagreements on the need and use of infinitesimals and mathematicians began to worry about the lack of rigorous foundations of the calculus (recall that the foundations of mathematics is a sound reasoning about the infinitesimal).
              This would change due 
                                  Cauchy's  ideas of \emph{function and limit}  which led to a more  \emph{rigorous} formulation of the calculus, limit/continuity/\emph{real numbers}. And, due to the emerging \emph{exact definition of real numbers}  the rules for reasoning with real numbers became even more precise.   However, all this historical background of the development of analysis is rarely reflected in the modern teaching of the subject.  Instead, students are introduced to methods that they find challenging, like  
                                  the  $\epsilon-\delta$/ $\epsilon-N$ proofs of limits  without background material on why limits, infinites, approximations  and infinitesimals were developed.    From our experience, an evolutionary and somewhat historic approach is helpful.  This is why we embarked on a book~\cite{Kamareddine-seldin-analysis} that introduces mathematical analysis  by  employing the evolution of this area of mathematics to first develop fundamental concepts of mathematical analysis and to only introduce formal definitions after the concepts are understood.

                                  The landscape of mathematics would change forever during the 19th century and the commitment to rigorous foundations would lead to the discovery of computability and its limits.  Rigorous foundations also 
                                  shed light on the holes that                                   started to appear in Euclid's historic work which led some to question the deductive structure of the Elements.
   Such logical inaccuracies have been addressed in the work of Hilbert~\cite{Hilbert-geometry}  who wrote 20 postulates adequate to prove all the theorems in the Elements. Here we go through some of these holes.

  \begin{itemize}
   \item
     Look at Proposition 1 of Book I of Euclid’s Elements:
\begin{tcolorbox}[width=.88\textwidth, colframe=red]
\emph{To construct an equilateral triangle on a given finite straight line.}
\end{tcolorbox}
Let $AB$ be the finite straight line. The proof draws two circles
with radius $AB$, and center $A$ (resp.\ $B$). The circles intersect at  $C$ and the triangle $ABC$ is equilateral. 
\begin{center}
 \begin{tikzpicture}[scale=0.3]
   \path (0,2) coordinate (A) (4,2) coordinate (B);
     \draw (0,2) circle (4)node[left]{A};
     \draw (4,2) circle (4)node[right] {B};
     \node at (-4.5,2) {D};
     \node at (8.5,2) {E};
     \draw(A)--(B);
        \node at (2,6) {C};
 \end{tikzpicture}
 \end{center}
There is a problem in this proof.  At first glance, there does not appear to be
any doubt that the construction given there constructs the desired equilateral triangle
and that the proof proves that it is an equilateral triangle. However, there is a gap
in the proof. There is, in fact, no proof that the point $C$ exists.
We can  construct a \emph{model} of geometry in which all of the postulates and axioms are satisfied but Proposition 1 is not.
\item
  Euclid’s Postulate 5 (the parallel postulates) is less obvious than the other postulates. 
\item
  Euclid used a number of statements  as facts in his Elements even though they had neither been proved nor been introduced as postulates. For example:
    \begin{tcolorbox}[width=.88\textwidth, colframe=red]
         A straight line that intersects one side of a triangle but does not pass through any vertex of the triangle must intersect one and only one of the other sides.
    \end{tcolorbox}

  Based on this statement, Pasch proved that Euclid’s formulation was not complete in the sense that there are statements that should hold but which cannot be proven from Euclid’s formulation.
   \begin{enumerate}
  \item
    A straight line passing through the center of a circle must intersect the circle.
    \item
    Given 3 different points on the same line, one of them is between the other two.  
   \end{enumerate}
\end{itemize}

  Having introduced the discrete (natural, rational and integer) numbers, and having emphasised the historical treatment of  continuous magnitudes and the need for real numbers in the development of analysis,  we now discuss the real numbers.
  \subsection{What are the real numbers?}
  Recall Proposition 2.\ of Book VII of the \emph{Elements} and  the approximations for $\sqrt{2}$ in Section~\ref{numbmagsec}. You can think of $\sqrt{2}$ as all the rational numbers strictly less than it. I.e.,  as: $\{1, 1+ \myfrac{1}{2}, 1+\myfrac{1}{2+\myfrac{1}{2}},
  1+\myfrac{1}{2+\myfrac{1}{2+\myfrac{1}{2}}}, \cdots\}$.  All irrational numbers have infinitely distinct approximations like $\sqrt{2}$.  Hence, the real numbers can be defined as non empty subsets of the rationals which satisfy some properties (see below).    Real numbers will be defined as elements of a complete ordered field.  Hence the following 
    definitions.
  \begin{definition}
    A   \emph{field} is a set $(S, +, \cdot)$ such that 
      $S$ is closed under $+$ and $\cdot$ and satisfies distributivity $a(b+c) = ab+ac$,
      commutativity and 
      associativity of $+$ and $\cdot$,  and existence of identity elements $0$ and $1$ ($a+0= a$ and $a\cdot 1 = a$) and inverses $-a$ and $a^{-1}$ (for each  $a$ except for $0$ under $\cdot$).
  \end{definition}
  \begin{example}
    None of $\mathbb{N}^+$ or $\mathbb{Z}$ is a field but  $\mathbb{Q}$ is a field.
    \end{example}
  \begin{definition}
     A field is \emph{ordered} (by $<$) if for all  $a$, $b$, $c$:
            \begin{itemize}
\item exactly one of $a < b$, $a = b$, and $b < a$ holds.

\item if $a < b$ and $b < c$, then $a < c$.

\item if $0 < a$ and $0 < b$, then $0 < a + b$ and $0 < ab$.

\item $a < b$ if and only if $0 < b +(- a)$.
\end{itemize}
    \end{definition}
 
  The next axiom 
  is important for the real numbers.
  \begin{center}   \begin{tcolorbox}[width=.88\textwidth, colframe=red]
     \textsc{Axiom of Completeness [AC]}\\
Every nonempty set of quantities that has an upper bound has a least upper bound.
\end{tcolorbox}\end{center}
Now we give the definition of the {\bf Real Numbers $\mathbb{R}$}.
      \begin{definition}
    Our quantities form an ordered field that satisfies the Axiom of Completeness AC.  We will refer to them as {\it real numbers} and denote their collection  by $\mathbb{R}$.
      \end{definition}
      Recall that the real numbers are continuous whereas the natural/integer/rational numbers are discrete.  The following help us to see some differences between these numbers.
      \begin{center}   \begin{tcolorbox}[width=.88\textwidth, colframe=red]
           \textsc{Archimedes Law} [AL]\\  For any two quantities $a$ and $b$ where $b > a>0$, there is a positive integer $n$ such that $b < an$.
      \end{tcolorbox}\end{center}
      \begin{definition}
        An ordered field which also satisfies AL is called an \emph{Archimedean ordered field}.
        \end{definition}
      \begin{example}
        $\mathbb{Q}$ is an Archimedean ordered field.
        \end{example}
      \begin{itemize}
      \item
  {\bf Completeness implies the Archimedean Property}
  Assume $a$ and $b$ are real numbers such that  $a > 0$.
  There is a positive integer $n$ such that $an > b$.
\item
   We can approximate real numbers by rational numbers.

{\bf Density of rationals}
If $a$ and $b$ are any two real numbers with $a < b$, then there is a rational number $r$ such that $a < r < b$.
\end{itemize}


\begin{thebibliography}{00}
  \bibitem{berkeley:Anls}
George Berkeley.
\newblock {T}he {A}nalyst: or {A} {D}iscourse {A}ddressed to an {I}nfidel
  {M}athematician.    {F}irst
  printed in 1734.
\newblock In {A. A.} Luce and {T. E.} Jessop, editors, {\em The Works of
  {G}eorge {B}erkeley {B}ishop of {C}loyne}, volume~4, pages 53--102. Nelson,
  London, 1951.
\bibitem{Heath:BEE-1956}
Heath.
\newblock {\em The 13 Books of Euclid's Elements}.
\newblock Dover, 1956.

\bibitem{Hilbert-geometry}
David Hilbert.
\newblock {\em The Foundations of Geometry}.
\newblock The Open Court Publishing Co, 1902.
\bibitem{Kamareddine-seldin-analysis}
  Fairouz Kamareddine and Jonathan Seldin.
  \newblock {\em A Primer of Mathematical Analysis and the Foundations of Computation.}
  College publications, ISBN 978-1-84890-443-9, October 2023. 434 pages.
\bibitem{knorr:EEEl}
{W. R.} Knorr.
\newblock {\em The Evolution of the Euclidean Elements: A Study of the Theory
  of Incommensurable Magnitudes and Its Significance for Early Greek Geometry}.
\newblock Reidel, Dordrecht and Boston and London, 1975.
\bibitem{rashed:arabmaths}
  Roshdi Rashed.
  \newblock {\em The development of Arabic mathematics: between arithmetic and algebra}.
  \newblock    Boston Studies in the Philosophy and History of Science (BSPS, volume 156).
  1994.
  
\bibitem{rashed:arabmathsfrench}  
  Roshdi Rashed.
  \newblock {\em Entre arithmétique et algèbre: Recherches sur l'histoire des mathématiques arabes}.
  \newblock Ouvrage publié avec le concours de l'Unesco. Société d'édition "Les Belles Lettres".
  Paris, 1984.
\end{thebibliography}
\end{document}